\documentclass[12pt,a4paper]{article}
\usepackage{amsmath,amssymb,amscd, float, verbatim, latexsym,subfig,cite,enumerate,setspace,mathtools}
\usepackage[vmargin=1in,hmargin=1in]{geometry}
\usepackage{amsthm,hyperref,amstext}
\usepackage[all]{xy}
\usepackage{lscape,ulem}
\usepackage{color} 
\usepackage{hyperref}

\newtheorem{theorem}{Theorem}[section]

\newtheorem{conjecture}[theorem]{Conjecture}

\newcommand{\R}{\mathbb{R}}
\newcommand{\C}{\mathbb{C}}

\newcommand{\N}{\mathbb{N}}

\newcommand{\Z}{\mathbb{Z}}
\newcommand{\T}{\mathbb{T}}

\newcommand{\cM}{\mathcal{M}}

\newcommand{\cW}{\mathcal{W}}

\title{Mechanisms of unstable blowup in a quadratic nonlinear Schr\"odinger equation}

\author{	Jonathan~Jaquette \footnote{Department of Mathematical Sciences, New Jersey Institute of Technology; Newark, NJ, USA; jcj@njit.edu} \footnote{Department of Mathematics and Statistics, Boston University; Boston, MA, USA}
}
\begin{document}

\maketitle

\begin{abstract} 
	In  the work 	
	Cho et al. [\emph{Jpn. J. Ind. Appl. Math.} 33 (2016): 145-166] the authors conjecture that the quadratic  nonlinear Schr\"odinger equation (NLS) $	i u_t =    u_{xx} + u^2 $ for $ x \in \mathbb{T}$ is globally well-posed for real initial data. 
We identify initial data whose numerical solution blows up in contradiction of this conjecture. 
The solution exhibits self-similar blowup and  potentially nontrivial self-similar dynamics, however the proper scaling ansatz remains elusive.

	Furthermore, the set of real initial data which blows up under the NLS dynamics appears to occur on a codimension-1 manifold, and we conjecture that it is precisely the stable manifold of the zero equilibrium for the nonlinear heat equation $	 u_t =  u_{xx} + u^2 $. 
	We apply the parameterization method to study the internal dynamics of this manifold, offering a heuristic argument in support of our conjecture. 
	 
\end{abstract}

\bigskip

%{\bf Keywords : } Nonlinear Schr\"odinger equations, integrable system, (quasi)periodic orbits, finite time blowup,  nonconservative equation.

%\bigskip
%\bigskip
%\centerline{{\bf AMS subject classifications}}
%\medskip
%\centerline{
%	35B10,  %	Periodic solutions to PDEs
%	35B44, %  	Blow-up in context of PDEs
%	35Q55,  	%   NLS equations (nonlinear Schrödinger equations) {For dynamical systems and ergodic theory, see 37K10}
%	37K10 % Completely integrable infinite-dimensional Hamiltonian and Lagrangian systems, integration methods, integrability tests, integrable hierarchies (KdV, KP, Toda, etc.)
%
%}  

%TODO List
%
%\begin{itemize}
%	  \setlength{\parskip}{0pt}
%	\setlength{\itemsep}{1pt plus 1pt}
%	\item Read \cite{fiedler2024real}
%%	\item Remake Figure \ref{fig:ContinuationInZeroAverage} (Continuation in zero average)
%	\item Remake Figure \ref{fig:BlowupPlot} (Blowup plots)
%	\item Clean up Section \ref{sec:Parameterization} (Parameterization) 
%%	\item Clean up the code for computing the parameterization 
%%	\item Remake Figure \ref{fig:ToySystemSimulation} (Toy blowup) 
%\end{itemize}

\section{Introduction}

In this paper we present   evidence to the contrary of a conjecture made by Cho et al.: 
\begin{conjecture}[\cite{Cho2016}]	\label{conj:GWP}
	The nonlinear Schr\"odinger equation (NLS):
	\begin{align} \label{eq:NLS}
		i u_t &=   u_{xx} + u^2 
	\end{align}
	is globally well-posed for any real initial
	data, small or large.
\end{conjecture}
 
The motivation for this conjecture arose in their work  continuing solutions  past the point of blowup in the following  nonlinear heat equation: 
\begin{align} \label{eq:Heat}
	u_t &= u_{xx} + u^2, & x \in \T .
\end{align}
Equation \eqref{eq:Heat} defines an analytic semi-flow,  hence solutions may be analytically continued for complex values of time \cite{sell2002dynamics}. 
If some initial data $u(0)=u_0$ blows up at finite  time $T$, one may attempt to bypass it through a detour in the upper/lower complex half-plane.  
The  solution to \eqref{eq:Heat} along the contour $ \Gamma^\theta = \{t = \tau   e^{i \theta}   \in \C:  0 \leq \tau  \}$ may be obtained by solving the following PDE:
\begin{align} \label{eq:Family}
	u_\tau &= e^{i \theta} \left(u_{xx} + u^2\right) ,&x \in \T .
\end{align}
so long as $ \theta \in [-\pi/2,\pi/2]$. 
Several such contours may be concatenated.  If no singularity occurs,  then the solution $u(T)$ at complex time $T\in\C$ will not depend on the path taken between $ 0$ and $T$. 

In earlier work by Masuda,  he showed that if small non-constant initial data is continued around the blowup point it will result in a branching singularity \cite{MASUDA1983119,Masuda1984}. For example, if a solution  $u(t)$ blows up at real time $T$ then  one can define the solution $u^{\pm}(2T)$   via two routes -- by following a path in either the upper or lower complex half plane -- and these two definitions will be  distinct $u^{+}(2T) \neq u^{-}(2T)$.  
Cho et al.  numerically investigated this phenomenon for large initial data and similarly observed branching singularities  \cite{Cho2016}, and there has been a growing interest in studying complex time singularities in PDEs  \cite{stuke2018complex,takayasu2022rigorous,fasondini2024complex,fiedler2024real,fiedler2024scalar}. 

Cho et al. further observed that all the real initial data they simulated under \eqref{eq:NLS} would decay to zero, hence their Conjecture \ref{conj:GWP}. 
%Many of these individual numerical results were subsequently confirmed  in \cite{takayasu2022rigorous} through the use of computer-assisted-proofs, and 
Supporting evidence toward this conjecture may be found \cite{jaquette2022global} where it was shown that  close-to-constant real initial data will  decay to zero in both forward and backward time. 
%\begin{theorem}[Theorem 1.3 \cite{jaquette2022global}]
%	\label{thm:TrappingRegion}
%	Fix a real scalar $y_0  \in \R$ and a function $ w_0 :\T^1 \to \C$ on a $ 2\pi/\omega$ periodic torus having summable Fourier coefficients, that is $ w_0 (x)= \sum_{k \in \Z} a_k e^{i k \omega x}$ for $ a = \{a_k\} \in \ell^1$. Let $ u(t)$ be the solution of \eqref{eq:NLS} with initial data $u_0(x) = y_0 +w_0(x)$. If 
%	\(
%	\| a\|_{\ell^1} < e^{- \pi/2} |y_0|, 
%	\) 
%	then $ \lim_{t \to \pm \infty} u(t) = 0$.  
%\end{theorem}
\begin{theorem}[Theorem 1.3 \cite{jaquette2022global}]
	\label{thm:TrappingRegion}
	Fix a complex scalar $z_0 =r e^{i \phi} \in \C$ and a function $ w_0 :\T^1 \to \C$ on a $ 2\pi/\omega$ periodic torus having summable Fourier coefficients, that is $ w_0 (x)= \sum_{k \in \Z} a_k e^{i k \omega x}$ for $ a = \{a_k\} \in \ell^1$. Let $ u(t)$ be the solution of \eqref{eq:NLS} with initial data $u_0 = z_0 +w_0$, and suppose: 
	\[
	\| a\|_{\ell^1} < e^{- \pi/2} |z_0|.
	\] 
	If $0 \leq \phi \leq \pi$ then $ \lim_{t \to - \infty} u(t) = 0$. 
	If $\pi  \leq \phi \leq 2\pi$ then $ \lim_{t \to + \infty} u(t) = 0$. 
\end{theorem}

Recent work by Fiedler and Stuke \cite{fiedler2024real}
may be seen to provide contradictory evidence to Conjecture \ref{conj:GWP}. 
In this work they study  the family of PDEs:
\begin{align} \label{eq:Fiedler}
	u_t &=  e^{i \theta}   
	\left(
	u_{xx} + 6u^2  - \lambda  \right) ,
	&x\in \T
\end{align}
and in essence show that real eternal solutions are not complex entire.  They prove there exists real initial data to \eqref{eq:Fiedler} which for the heat equation $ (\theta=0)$ is a heteroclinic orbit existing globally in time, but for the NLS $(\theta = \pm \tfrac{\pi}{2}$) it blows up in finite time. 
However, this analysis requires $\lambda >0$ 
and does not resolve Conjecture \ref{conj:GWP}.

While \eqref{eq:NLS} is  formally  a nonlinear Schr\"odinger equation, it  is quite different from the typical focusing/defocusing NLS   of the form $ i u_t  = \triangle u \pm |u|^{p-1} u$. Most notable is the lack of gauge symmetry or any disernable Hamiltonian form.  Moreover  \eqref{eq:NLS} does not preserve any nontrivial analytic conserved quantities  \cite{jaquette2022global}. 
Such nonstandard NLS have primarily garnered attention in works studying the minimal regularity needed to establish local well-posedness  \cite{kenig1996quadratic,bejenaru2006sharp,iwabuchi2015ill}. 
The particular form of the nonlinearity has a large impact on the well-posedness theory; the equation $  i u_t +\triangle u =  |u|^p$, similarly lacking gauge invariance, has been shown to  blowup for  a large class of initial data 
\cite{oh2012blowup,ikeda2015some,fujiwara2017lifespan,fujiwara2021global}.

As a complex PDE, the nonlinearity in   \eqref{eq:NLS} and \eqref{eq:Heat} is closely related to the 1-dimensional Constantin-Lax-Majda  model   of vortex stretching in an incompressible fluid: 
\[
\omega_t = \omega H(\omega) 
\]
where $H$ is the  Hilbert transform \cite{constantin1985simple}. 
The Hilbert transform is a nonlocal, singular integral transform and, importantly for our purposes, $H^2 (\omega)=-\omega$. 
By making the complex change of variables 
$u(t,x)=H(\omega(t,x)) + i \omega(t,x)$, one obtains the local equation: 
\begin{align}
	u_t =\frac{1}{2} u^2 \label{eq:CLM_Local}
\end{align}
which is 
uncoupled from all spatial dependence. 
This equation can readily be seen to produce blowup in finite time, however it is missing key features of a model for incompressible fluids which may serve to prevent blowup. 
In \cite{de1990one,de1996partial} de Gregorio further incorporated terms to account for nonlinear advection and viscosity yielding the eponymous equation, which  has since been   modified and widely studied \cite{okamoto2008generalization,hou2008dynamic,guo2013convergence}. 

\begin{figure}[t]
	\begin{center}
		\includegraphics[width=1\textwidth]{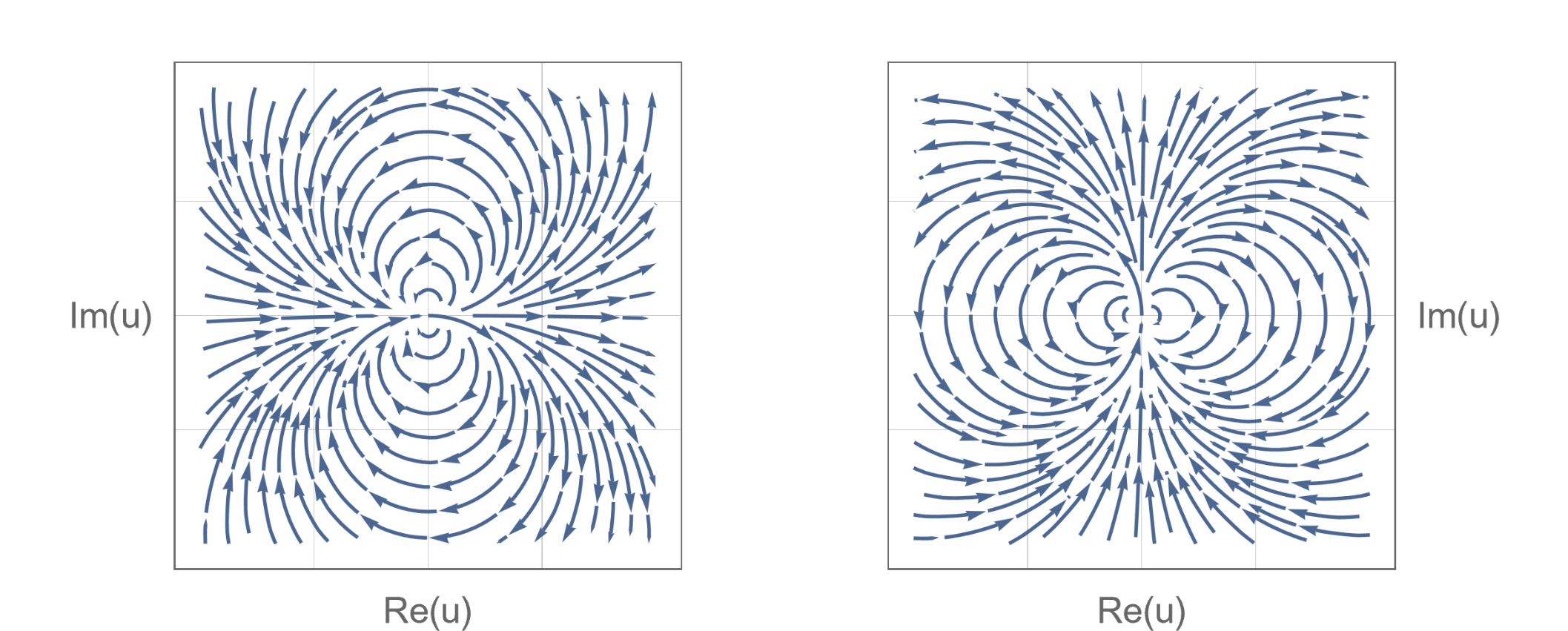}
	\end{center}
	\caption{ The complex phase portraits for ODEs $ u_t = u^2 $ (left) and $i u_t = u^2 $ (right).}
	\label{fig:StreamPlots}
\end{figure}

As a real-valued equation \eqref{eq:CLM_Local} is a textbook example of a differential equation leading to blowup for   positive initial data. However richer dynamics may be found for complex initial data, see Figure \ref{fig:StreamPlots}. 
As \eqref{eq:CLM_Local} has only a nominal spatial dependence, its solution for initial condition $ u(0,x)=u_0(x)$ may be found explicitly:
\begin{align*}
	u(t,x)&= \frac{u_0(x)}{ 1 - \tfrac{1}{2} t u_0(x)}
\end{align*} 
This solution will blow up in forward (backwards) time if there is some $x$ for which $u_0(x)$ is a positive (negative) real number. Conversely if $ u_0(x) \in \C \backslash \R \cup \{ 0\}$ then it is globally well-posed. 
Such a case where blowup only occurs on a meagre set is difficult to observe through direct numerical simulation, and the development of techniques to find elusive blowup solutions is an active area of research \cite{protas2022systematic,wang2023asymptotic}. 
On the other side of the coin one can view the problem of blowup probabilistically, and answer questions of local and global well-posedness with respect to an associated
Gibbs measure \cite{bourgain1996invariant,liu2023probabilistic,deng2019invariant}.

Returning to Conjecture \ref{conj:GWP}, one route to proving global well-posedness could be achieved by showing that all real initial data limits to zero. 
While this might be the case for almost all initial data, the NLS in \eqref{eq:NLS} exhibits a rich dynamical structure beyond decay to zero \cite{jaquette2022global,feng2023long,jaquette2022quasiperiodicity,jaquette2022singularities}.   
There are two families of nontrivial equilibria to \eqref{eq:NLS} given by the Weierstrass elliptic functions \cite{fiedler2024real,jaquette2022global}, with heteroclinic orbits connecting the trivial and nontrivial equilibria \cite{jaquette2022global,jaquette2022singularities}. 
In \cite{jaquette2022quasiperiodicity} it is shown that small initial data supported on positive Fourier modes will yield periodic orbits of locked frequency. 
For larger initial data, however, global existence cannot be assured: 
\begin{theorem}[Theorem 1.5\cite{jaquette2022quasiperiodicity}] \label{prop:Intro_Blowup}
	Consider   \eqref{eq:NLS}  with   	   initial data $u_0(x) =   \alpha  e^{ i \omega  x}$ with  $\alpha \in \C$,  $ \omega >0$.   
	\begin{enumerate}[(a)]
		\item If $|\alpha| \leq  3  \omega^2$, then  $u$   is a smooth periodic solution with period $\frac{2 \pi}{\omega^2}$. 
		\item If $|\alpha| \geq  6\omega^2$,  then $\|u(t)\|_{L^2}$  blows up  in finite time  $|T_*| \leq  \frac{2 \pi}{\omega^2}$.
		%		In particular, there exists some  $|T_*| \leq  \frac{2 \pi}{\omega^2}$ such that  $\limsup_{t \to T_*} \| u(t)\|_{L^2} = + \infty$.   
	\end{enumerate}
\end{theorem}
%To compare with Theorem \ref{thm:TrappingRegion}; there solutions that are close to constant / zero average will decay to zero in forward and backwards time. In contrast, 
Theorem \ref{prop:Intro_Blowup} presents a dichotomy of solutions, between those which exist globally and are periodic if the initial data is small, and those which will blowup in finite time if the initial data is large. 
While this theorem helps elucidate the geometry of solutions which blowup, the initial data is complex valued and does not directly help with resolving Conjecture \ref{conj:GWP}.

In Section \ref{sec:EvidenceOfBlowup} we search for real initial data leading to blowup by analyzing a 1-parameter family of  initial data: $u_0(x;A) = 30 \cos(2 \pi x) + A$.  The blowup result in Theorem \ref{prop:Intro_Blowup}, occurring for a function with zero average, hints that small values of $A$ may lead to blowup. More concretely, by Theorem \ref{thm:TrappingRegion} 
	 only values $A \in  [-150,150]$ need to be considered, as larger values of $|A|$ are guaranteed to converge to zero. Indeed, most of the solutions we simulate in this range converge to zero, however for a fine enough sampling of parameters  we can observe solutions which do not decay. 
All source code is available at   \cite{codesUnstableBlowupNLS}.

Further parameter tuning reveal a solution which, at least numerically, blows up in finite time. 
Unlike the blowup in Theorem \ref{prop:Intro_Blowup} which occurs in time $ |T_*| <1/2 \pi$, the blowup solution slowly grows in size over a long time before blowing up.  Throughout this process the solution oscillate with increasing vigor and complexity, pumping energy into the higher Fourier modes.

Surprisingly if we take these same initial conditions which blowup under the NLS dynamics and numerically integrate them under the heat equation, they appear to  approach the zero equilibrium exponentially fast. Such behavior is not generic under the dynamics of \eqref{eq:Heat}.
If $A >0$, then the solution is guaranteed to blowup in finite time. 
 Whereas if $A \ll 0$, then the solution will be exponentially attracted to the center manifold, along which it will converge algebraically to zero \cite{takayasu2022rigorous}. 
The strong stable manifold $ \cW^s(0)$ of \eqref{eq:Heat} is codimension-1 and divides these two types of behavior. Hence we have motivation for the following conjecture: 

\begin{conjecture}	\label{conj:ParabolicManifold} 
	Let $ \cW^s(0) \subseteq C(\T^1,\R)$ denote the strong stable manifold of the zero equilibrium for the nonlinear heat equation \eqref{eq:Heat}. 
	Initial data $ u_0 \in C(\T^1,\R)$ with summable Fourier coefficients is globally well-posed under \eqref{eq:NLS} if and only if $u_0 \notin \cW^s(0) $. 
\end{conjecture}

In Section \ref{sec:Resonances} we present a heuristic argument in support of this conjecture by analyzing the internal dynamics on the manifold $ \cW^s(0)$ for finite Galerkin truncations of the NLS equation.  
In order to compute the strong stable manifold we use the parameterization method to obtain a Taylor series approximation of both the manifold and its internal dynamics \cite{cabre2003parameterizationI,cabre2003parameterizationII,cabre2005parameterizationIII}.  
Since equations \eqref{eq:NLS} and \eqref{eq:Heat} only differ by multiplication by a complex scalar, our approximation of the strong stable manifold in the nonlinear heat equation immediately yields an approximation for an invariant   manifold of the nonlinear Schr\"odinger equation, of course with different internal dynamics. 

In general, if the stable eigenvalues of an equilibrium are rationally independent (i.e. non-resonant),  then the internal dynamics of the stable manifold are smoothly conjugate to the linearized  dynamics.  
The zero equilibrium in \eqref{eq:Heat} essentially  has a linearization with eigenvalues $-n^2$ for $ n \in \Z$, and resonances abound.  
This does not greatly affect the internal dynamics of the stable manifold to \eqref{eq:Heat}.  However in the center manifold of \eqref{eq:NLS} these resonances give rise to secular growth in the higher modes, explaining the chaotic cascade of energy to the higher modes observed in Section \ref{sec:EvidenceOfBlowup}. 
This analysis may be seen to be in analogy with works such as  \cite{colliander2010transfer,guardia2015growth} which,  by analyzing a resonant system of (infinitely many) ODEs, shows that there exists arbitrarily small initial data to the defocusing NLS on $\T^2$ which will exhibit arbitrarily large finite growth in the higher Sobolev norms. 
While our analysis offers a heuristic explanation for how secular growth  causes small initial data on $\cW^s$ to become large, we are unable to close this argument to provide a proof of blowup.

To robustly  analyze blowup solutions, it is common to make a self-similar ansatz for  scaling the magnitude and spatial dependence of a solution as it approach the  blowup time \cite{eggers2008role,quittner2019superlinear}.  
This self-similar change of variables renormalizes the blowup problem into one of studying the self-similar dynamics of a new PDE produced by the change of variables. 
In the simplest case, a blowup profile can be equated with an  equilibrium of a new PDE. 
However the self-similar dynamics need not be simple, and may in fact be periodic (also referred to as discretely self-similar) \cite{martin2003global,tao2016finite}, or even chaotic \cite{eggers2008role,campolina2018chaotic}.

In Section \ref{sec:SelfSimilar} we use self-similar variables to analyze the blowup of solutions to \eqref{eq:NLS}. 
The blowup solution observed in Section \ref{sec:EvidenceOfBlowup} exhibits complex spatio-temporal oscillations, suggestive of nontrivial self-similar dynamics, however  we shall side-step this complicated behavior.  
Instead, we study the solution starting from the larger initial data $ u_{300}(x) = 300 \cos(2 \pi x) + A_{300}$. Informed by Conjecture \ref{conj:ParabolicManifold}, we select the constant $ A_{300}$ so that $u_{300}$ is in the strong stable manifold $ \cW^s(0)$ of the heat equation. As predicted, this choice yields real initial data which numerically blows up under the dynamics of \eqref{eq:NLS}. 
By starting with larger initial data the solution blows up more directly, without as long of a period of secular drift exciting the higher Fourier modes. 
For  comparison we also study the  monochromatic initial data $u^{mc}_{300}(x) = 300 e^{2 \pi i x}$, which is guaranteed  to blowup by Theorem \ref{prop:Intro_Blowup}. 

Both solutions, starting from real initial data and from monochromatic initial data, appear to blowup in a self-similar manner. 
While the two blowup profiles look similar, they seem to obey different scaling rates.  
For the real initial data, the solution first   approaches a blowup profile with a single stationary blowup point, and later limits to a profile with two, non-stationary blowup points. 
It remains unclear what the proper self-similar ansatz needs to be. 
Furthermore,  even if the existence of such a self-similar blowup were to be established, for Conjecture \ref{conj:GWP} to be resolved it would still require showing that the stable set of such a blowup profile includes real initial data. 
\\

%%%%%%%%%%%%%%%%%%%%%%%%%%
%%%%%%%%%%%%%%%%%%%%%%%%%%

%\clearpage

\section{Evidence of blowup from real initial data } 
\label{sec:EvidenceOfBlowup}

\begin{figure}[b!]
	\begin{center}
		\includegraphics[width=1\textwidth]{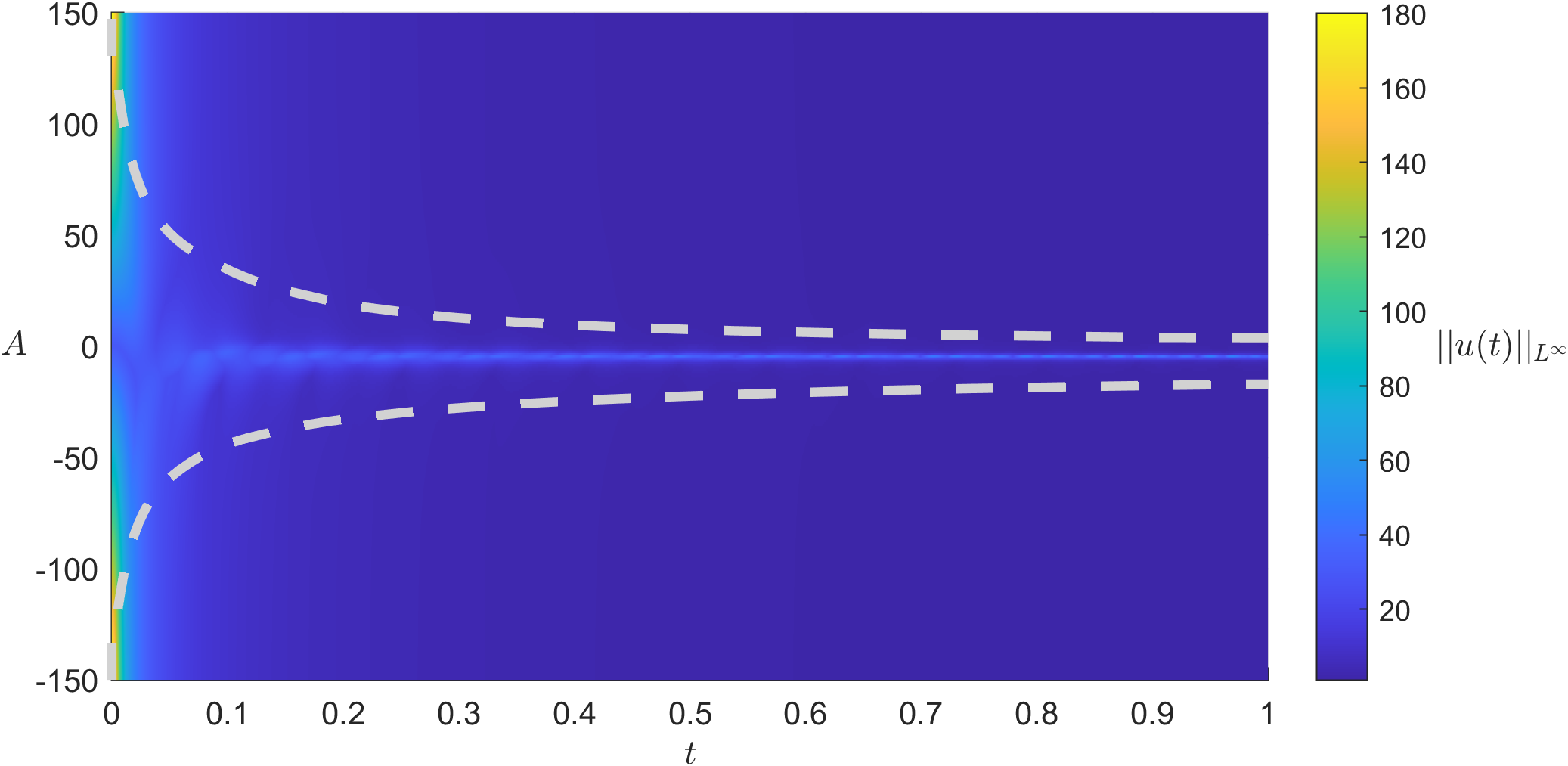}
	\end{center}
	\caption{ The norm $\|u_{30}(t;A)\|_{L^\infty}$ for  $ u_{30}(x;A)$  as in \eqref{eq:InitialDataFamily}. The dashed gray lines represent where the solution reaches the trapping region given by Theorem \ref{thm:TrappingRegion}.  }
	\label{fig:ContinuationInZeroAverage}
\end{figure}

Numerical evidence suggests that most real initial data to \eqref{eq:NLS} is globally well posed.   Indeed, it follows from Theorem \ref{thm:TrappingRegion} that for real initial data which is  close-to-constant (and arbitrarily large!) will limit to zero in both forward and backwards time.  
In light of this, to look for possible counterexamples to Conjecture \ref{conj:GWP}  one must consider initial data which is close to zero average. 
To that end, consider the 1-parameter family of real initial data: 
\begin{align} \label{eq:InitialDataFamily}
	u_{30}(x;A) &=  30 \cos( \tfrac{2\pi}{\omega } x) +A
\end{align}
Note that due to symmetry the Fourier series of the solution can be   expressed as a cosine series with complex, time-varying coefficients: 
\begin{align} \label{eq:CosineAnsatz}
	u(t,x) = \sum_{n \in \Z} a_n(t) e^{\frac{2 \pi i n}{\omega} x} =  a_0(t) + 2 \sum_{n=1}^{\infty} a_n(t) \cos(\tfrac{2 \pi n  }{\omega }x)
\end{align}
By Theorem \ref{thm:TrappingRegion} if $ |A| > 30 e^{\tfrac{\pi}{2}}  \approx 144.3$, then  the solution is guaranteed to exist for all $ t\in\R$ and converge to zero in both forward and backwards time. 
	In Figure \ref{fig:ContinuationInZeroAverage} is depicted the solution of the NLS  \eqref{eq:NLS} over the time interval $ t \in [0,1]$,  taking  initial data in \eqref{eq:InitialDataFamily}  with integers $ -150 \leq A \leq 150$.  This computation was performed using 256 Fourier modes and a time step of $h=10^{-4}$ with an exponential integrator \cite{driscoll2014chebfun}.

Almost all initial data converges to zero, however close to $A=-5$ the solutions take a long time to decay.  
For each solution which does decay to zero, Theorem \ref{thm:TrappingRegion} enables us to identify when it enters a so called ``trapping region'' where it is guaranteed to converge to zero. This is denoted in Figure \ref{fig:ContinuationInZeroAverage} by the dashed gray line.  
   Note that the time at which solutions enter the trapping region does not correlate with their $ L^\infty$ norm.  
 Indeed, the trapping region described by Theorem \ref{thm:TrappingRegion} is essentially a cone of Fourier series, whose zeroth Fourier mode both  has a negative imaginary part and is sufficiently larger than the higher Fourier modes.  

While a small neighborhood of initial data with $A$ about $-5$ does decay, it is more noticeable that nearby values of $A$ take asymptotically long to enter the trapping region.  
This provides us with a more robust measure of identifying solutions which do not converge to the zero equilibrium. 
We observe the same behavior for other families of initial data of the form $ u_g(x;A) = g(x) + A$, which suggests the existence of a codimension-1 manifold for which initial data to the NLS does not converge to zero.   
%Note, if we only look at the norm of solutions, then the spike at $ A=-5$ is  easy to miss. However by looking at the time it takes to reach the trapping region, we are better able to detect this. 

%This creates a sort of parabolic manifold, which contains zero as a equilibirum, which is stable (not asymptotically so) within the manifold. This seems to be a sort of Lyapunov family / folitation of periodic orbits.  
%Note here, this family of periodic orbits for the NLS are the stable manifold for the heat equation. 

%Also note that by only focusing on positive modes, we remove a resonance; there should really be a multiplicity two eigenvalue of number $k$. Compare this with what happens next; we are able to remove the multiplicity by making a cosine ansatz, but the resonance still remains. 

%%%%%%%%%%%%%%%%%%%%%%%%%%%%%%%%%%%%%%%%%%%%%
%%%%%%%%%%%%%%%%%%%%%%%%%%%%%%%%%%%%%%%%%%%%%
%%%%%%%%%%%%%%%%%%%%%%%%%%%%%%%%%%%%%%%%%%%%%
%%%%%%%%%%%%%%%%%%%%%%%%%%%%%%%%%%%%%%%%%%%%%
%%%%%%%%%%%%%%%%%%%%%%%%%%%%%%%%%%%%%%%%%%%%%
%\clearpage

\begin{figure}[b!]
	\begin{center} 
	\includegraphics[width=1\textwidth]{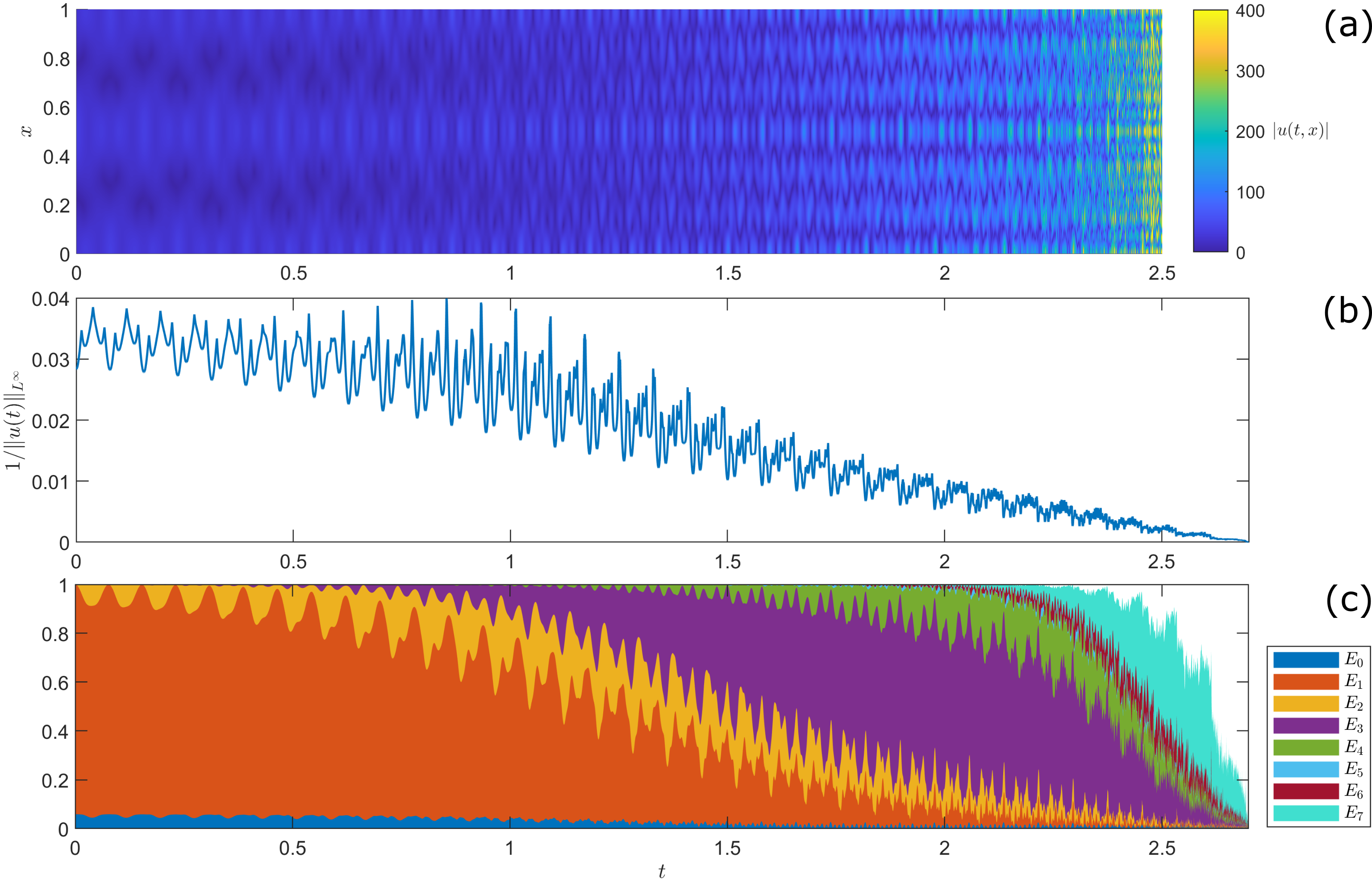}
	\end{center}
	\caption{ Solution to the NLS in \eqref{eq:NLS} with initial data \eqref{eq:InitialDataFamily} with 	$A_{30} = -5.3070235$: (a) Space time plot of the solution; (b)  Inverse norm of the solution $1/\| u(t)\|_{L^{\infty}}$;   (c) Relative proportions $E_n$ of the Fourier modes, see \eqref{eq:ProporionalEnergy}. }
	\label{fig:BlowupPlot}
\end{figure}

Informed by Conjecture \ref{conj:ParabolicManifold} we selected $ 	A_{30} = -5.3070235$ in \eqref{eq:InitialDataFamily}, and the solution of this initial value problem under the dynamics of the NLS \eqref{eq:NLS} is plotted in Figure \ref{fig:BlowupPlot}.  
This computation used a 4096 Fourier mode truncation and a time step of $h=10^{-7}$.
We selected $ A_{30}$ using a bisection method so  that $u_{30}$ would be on the stable manifold $\cW(0)$ for the nonlinear heat equation \eqref{eq:Heat}.   
%We note that the precise value of $A$ was selected with consideration made to Conjecture \ref{conj:ParabolicManifold}, so that $ u_{30}(x;A)$ would be contained in the stable manifold $\cW$ of $0$ under the dynamics of the nonlinear heat equation \eqref{eq:Heat}. 
We also note that nearby values of $A$ (e.g. $\approx\pm1 \%$) also appear to blowup under the dynamics of \eqref{eq:NLS}.

Overall, the trajectory appears to  oscillate while steadily growing larger.
As shown in Figure \ref{fig:BlowupPlot} (b), the $L^{\infty}$ norm of the solution  maintains regular oscillations of fixed period yet steadily grows larger in amplitude.  
In Figure \ref{fig:BlowupPlot} (c) we display a compositional breakdown of the $L^2$ norm according to the relative contribution from the first eight Fourier modes.   
By Parseval's identity  $ \|u(t)\|_{L^2}^2 = \|a(t)\|_{\ell^2}^2= |a_0(t)|^2 +2\sum_{n=1}^{\infty} |a_n(t)|^2$, and we plot in Figure \ref{fig:BlowupPlot}(c)  the ratio $E_n(t)$ of the $n^{\text{th}}$ mode: 
\begin{align} \label{eq:ProporionalEnergy}
	E_n(t):=
	\begin{dcases}
		 \frac{|a_n(t)|^2}{\|a(t)\|_{\ell^2}^2} & \mbox{if } n=0 \\
		 \frac{2|a_n(t)|^2}{\|a(t)\|_{\ell^2}^2} & \mbox{if } n\neq 0 
	\end{dcases}
\end{align}
Note that $ \sum_{n=0}^\infty E_n(t) =1$ by construction.

While the spacetime plot in Figure \ref{fig:BlowupPlot}(a) appears quite complex, the relative energy of each spatial mode in Figure \ref{fig:BlowupPlot}(c) paints a clearer picture. 
As the   solution evolves in time, higher spatial frequencies are progressively excited with increasing temporal oscillations. 
Furthermore,  these    oscillations grow ever more complex and almost fractal-like. 
%Overall, there appears to be two distinct qualitative features to this blowup solution: (i) on shorter time scales, the solution  periodic oscillates with progressive excitement of the higher modes; (ii) on longer time scales, the solution steadily grows, eventually leading to blowup.  

\section{Secular growth of solutions along a center manifold}
\label{sec:Resonances}

\subsection{Dynamics of the PDE's Galerkin truncation}

 To investigate the initial growth phase of the solutions we consider   finite Galerkin truncations of \eqref{eq:NLS}.   
% Note also, here we consider functions on a $2\pi$-periodic torus, whereas in other sections we consider functions on a $1$ periodic torus.  
For ease of analysis in this section we take $ \omega = 2 \pi$, and consider  a function given as a $2 \pi $-periodic  cosine series as in \eqref{eq:CosineAnsatz}.  
The $N$-mode truncation of  \eqref{eq:Family} yields the following dynamics on the Fourier modes: 
 	\begin{align} \label{eq:Galerkin}
 		e^{ -i \theta } \dot{a}_n &= -n^2 a_n +\sum_{\substack{n_1 + n_1 = n\\ -N \leq n_1, n_2 \leq N}} a_{|n_1|} a_{|n_2|}
 	\end{align} 
 	By the cosine ansatz we have $a_n = a_{-n}$, thereby this defines a complex ODE on $ \C^{N+1}$.  	
 Furthermore the system has  an equilibrium at $0\in \C^{N+1}$ whose linearization has eigenvalues $\lambda_n =  -e^{i \theta} n^2 $ for $ 0 \leq n \leq N$.  
 	
 	Define $ \cW(0)$ as the unique invariant manifold tangent to  $\{0\} \times \C^N$. 
 	 Selecting different values of $ \theta $ will induce different internal dynamics on $\cW$, however $\cW$ will remain an invariant manifold for any choice of $ \theta$.  For example,  if the real component of $-e^{i \theta}$ is negative, then $ \cW$ is a submanifold of the equilibrium's stable set, and trajectories on $\cW$ will exponentially approach the zero equilibrium  in forward time.
% 	  Conversely if $ \Re(e^{i \theta})<0$ then $ \cW$ is a submanifold of the equilibrium's unstable set. 

 	More generally if $\cM$ is any invariant complex manifold of any ODE $ \dot{x} =f(x)$, then $\cM$ will also be an invariant manifold of the system $ \dot{x}= \mu f(x)$ for scalars $ \mu \in \C$.  The most well know usage of this fact is the case $\mu = -1$. That is, the stable manifold of an equilibrium in forward time is also the unstable manifold of the equilibrium in backwards time. 
 	While trajectories on $\cM$ remain invariant when $ \mu \in \R$, this is not the case when $ \mu$ is complex. 
 	
% 	Thus one may consider the flow map $ x \mapsto \varphi(x, \mu t)$ as defining an automorphism on invariant manifolds $\cM$. 

 	\begin{figure}[t]
 		\begin{center}
% 			\hfill
 			\includegraphics[width=1\textwidth]{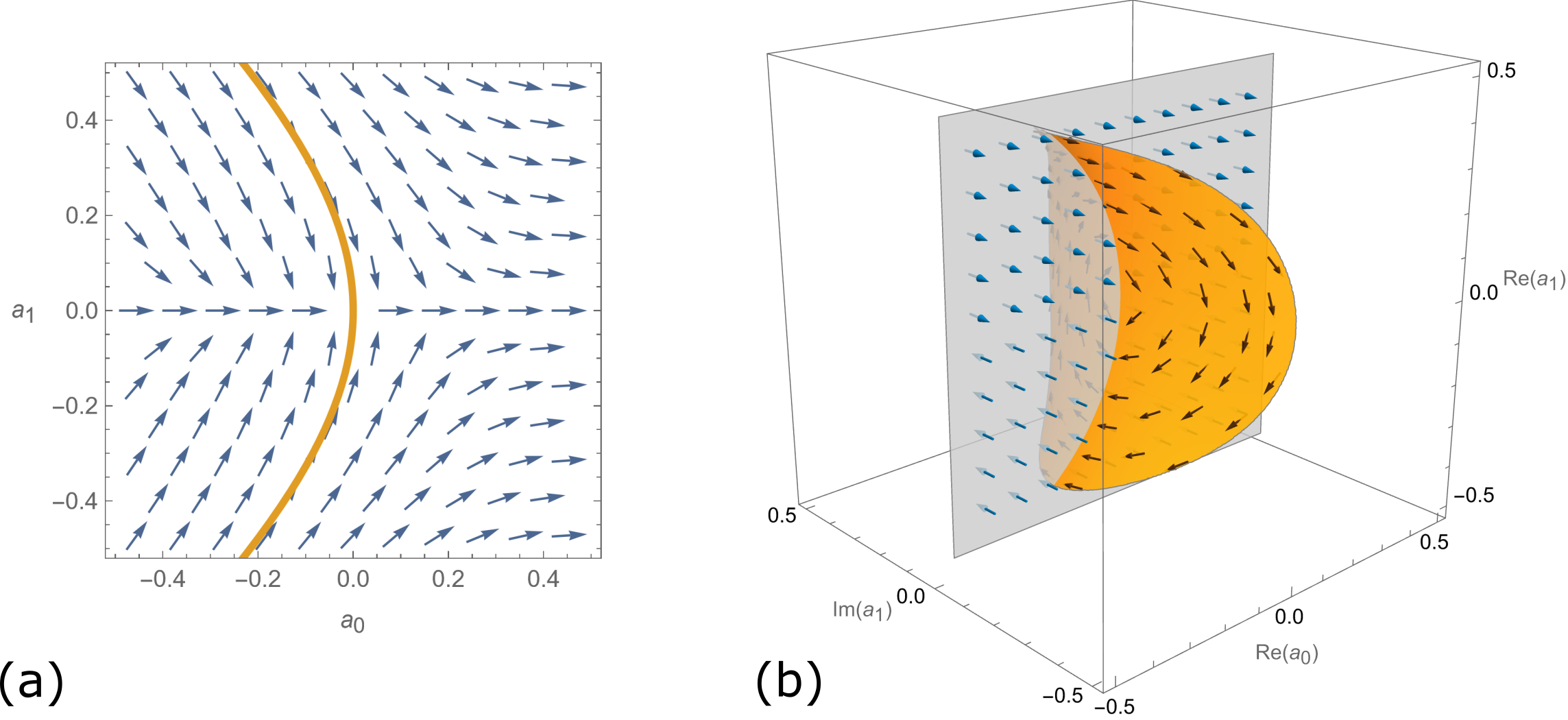}
% 			\includegraphics[width=.35\textwidth]{Figures/realWs}
% 			\hfill
% 			\includegraphics[width=.4\textwidth]{Figures/W2Mode_Imag}\hfill \null
 		\end{center}
 		\caption{Phase portraits of the dynamics to \eqref{eq:2Modes} with the invariant manifold $\cW$ tangent to $ \{0\}\times \C$ depicted in orange. (a) When $\theta=0$ real initial data remains invariant, and the stable manifold $\cW$ divides trajectories which converge to zero, from those which/that blowup in finite time. (b) When $\theta = \pi/2$, the invariant manifold $ \cW$ is foliated by periodic orbits.  }
 		\label{fig:GalerkinModel}
 	\end{figure}
 	
 	For a simple illustration,     consider the dynamics of \eqref{eq:Galerkin}    in the case $N=1$  where we obtain the system of complex ODEs: 
 \begin{align} \label{eq:2Modes}
 	e^{ -i \theta } \dot{a}_0 &= a_0^2 + 2 a_1^2 \\
 	e^{ -i \theta } \dot{a}_1 &= -a_1  + 2 a_0 a_1  \nonumber
 \end{align}
 
 The case $ \theta=0$ is analogous to the   nonlinear heat equation \eqref{eq:Heat},  see Figure \ref{fig:GalerkinModel}(a). 
 Real initial data is invariant, and if $ a_0>0$ then there is finite time blowup. 
  Furthermore the manifold $\cW$ divides solutions  which blowup in finite time from those which converge algebraically to zero. 
  Note also that  if $-\pi/2 < \theta < \pi/2$ then    solutions on $\cW$ will exponentially decay. 

In contrast, consider the case $ \theta=\pm \frac{\pi}{2}$  analogous to the NLS in \eqref{eq:NLS}, see Figure \ref{fig:GalerkinModel}(b). 
In this case the eigenvalues $ \{\lambda_0,\lambda_1\} = \{ 0, \mp i\}$ are purely imaginary and the manifold $\cW$ is a submanifold of the center manifold. 
Furthermore $ \cW(0)$ is foliated by a Lyapunov family of periodic orbits \cite{henrard1973lyapunov}. 
In the full PDE we can similarly observe periodic behavior for short time scales, see $ 0 < t < 0.25$ in Figure \ref{fig:BlowupPlot}. However, like with Figure \ref{fig:BlowupPlot}(c),  we will see in the dynamics of \eqref{eq:Galerkin} for larger values of $N$ that there is a  secular drift of solutions along $\cW$, whereby the lower   modes progressively excite the higher   modes.

\subsection{Parameterizing the manifold $\cW$  }
\label{sec:Parameterization}

The existence of stable, unstable, and center manifolds associated to equilibria has been established for a wide variety of dynamical systems, such as ODEs, PDEs, and DDEs \cite{sell2002dynamics}. 
However even if such a manifold is unique, its representation via a coordinate chart is not unique, nor is there just one method for computing an approximation for the manifold. 
%For a basic approximation, in a small neighborhood about the equilibrium one may approximate the unstable (stable) manifold by the unstable (stable) eigenspace, and then globalize the manifold by integrating forward (backwards) in time. While this may work for an ODE, PDEs may not be time reversible!  

To analyze and compute the invariant manifold $\cW$ associated with \eqref{eq:Galerkin} we shall use the parameterization method \cite{cabre2003parameterizationI,cabre2003parameterizationII,cabre2005parameterizationIII}, which has had great success analyzing the dynamics of  PDEs 
\cite{reinhardt2019fourier,jaquette2022global,jaquette2022singularities,gonzalez2022finite,jain2022compute,opreni2023high}.     
This approach may be seen to be in contrast to the graph approach, where one fixes the representation of the invariant manifold as graph over the particular eigenspace and the internal dynamics need to be solved for. Instead, the parameterization method fixes the internal dynamics and solves for the coordinate chart as a map into the entire phase space. 

To  briefly review the parameterization method  we follow the presentation in  \cite{haro2016parameterization}. 
	Consider $ F: \C^n \to \C^n$  and the  differential equation given by: 
	\[
	\dot{z} = e^{i\theta} F(z)
	\] 
	having an equilibrium $ z_*$, and suppose the linearization $ DF(z_*)$ has simple eigenvalues $ \lambda_1    \dots \lambda_n$.  Fix $ d\leq n$, and let $L \in \C^{n \times d}$ denote a matrix whose $d$-columns are eigenvectors of $ DF(z_*)$ with eigenvalues $ \lambda_1 \dots \lambda_d$, and let $ V \subseteq \C^n$ denote the $d$-dimensional subspace  spanned by the column vectors of $L$. 
	The goal is to  look for a parameterization $z = W(\sigma): \C^d \to \C^n$ of the invariant manifold $\mathcal{W} $ tangent to $V$ at $z_*$.  
	The internal dynamics on the manifold are described by a vector field $ \dot{\sigma} = e^{i \theta} f(\sigma)$ on $\C^d$ with $ f(0)=0$, and the invariance equation is given by:
\begin{align} \label{eq:Invariance}
	 e^{i\theta}	F(W(\sigma)) &= DW(\sigma)  e^{i\theta}f(\sigma)
\end{align}
Note that for the purpose of solving the invariance equation we are able to  cancel the complex phase $ e^{i\theta}$ completely! 
%Computing invariant manifolds of degenerate fixed points is often challenging \cite{baldoma2008one,baldoma2020invarianti,baldoma2020invariantii}, however we are able to avoid much of these  difficulties. 
Hence, while the internal dynamics of $ \cW$ will be different for different values of $ \theta$, the Taylor series parameterization of $ W$ we compute will be exactly the same. 

	The parameterization method seeks to write $W$ and $f$ as power series in $ \sigma = ( \sigma_1,   \dots ,\sigma_d)$: 
\begin{align}
		W(\sigma) &= z_* + \sum_{k\geq 1 }W_k(\sigma),
		&
	f(\sigma) = \sum_{k\geq1} f_k(\sigma)	,
\end{align}
where each $ W_k(\sigma)$ (respectively $f_k(\sigma)$) is a homogeneous polynomial in $\sigma$ of degree $k$ with coefficients in $\C^n$ (respectively with coefficients in $\C^d$), that is:
\begin{align}
	W_k(\sigma) &= \sum_{k_1 + \dots+ k_d = k} W_{(k_1\dots k_d)} \sigma_1^{k_1} \cdots \sigma_d^{k_d},
	&W_{(k_1\dots k_d)}&\in \C^n
\end{align}
%		where $ W_k$  (respectively $f_k$) is an $n$-dimensional (resp $d$-dimensional) vector. 
		
		To enforce the invariant manifold $\mathcal{W} $ being tangent to the eigenspace $V$ at the equilibrium we fix $ W_1(\sigma)=L\sigma$ and $ f_1 (\sigma) = \Lambda_L \sigma$,  where $\Lambda_L $ is the $d \times d$ diagonal matrix with diagonal entries $ \lambda_1 \dots \lambda_d$.  
		The higher order coefficients in the power series may be recursively obtained by matching the order-$k$ terms  in \eqref{eq:Invariance}, yielding the cohomological equations: 
	\begin{align}
		DF(z_*) W_k(\sigma) - DW_k(\sigma) \Lambda_L \sigma - L f_k(\sigma) = - E_k(\sigma)
		\label{eq:Cohomological}
	\end{align}
	where: 
	\[
	E_k(\sigma) = [F(W_{<k}(\sigma))]_k - [DW_{<k}(\sigma)f_{<k}(\sigma)]_k
	\]
	is the order-$k$ ``error'' term.  
	For each $k$, both sides of \eqref{eq:Cohomological} is a homogeneous polynomial in the variable $\sigma$ of degree $k$. Note however that   $E_k$ only depends on 
	  coefficients whose order is strictly less that $k$.   
	Thus through \eqref{eq:Cohomological} the coefficients $W_k$ may be solved recursively  for orders $ k \geq 2$, subject to the choice of coefficients $f_k$.

	The simplest possible internal  dynamics would be the linear dynamics $ \dot{\sigma}= e^{i \theta} \Lambda_L\sigma$ wherein $f_k =0 $ for $ k\geq 2$.  
	In this case, fixing notation $\vec{\lambda} = ( \lambda , \dots , \lambda_d)$, it follows that for each multi-index $ \vec{k} =(k_1, \dots ,k_d)$  equation  \eqref{eq:Cohomological} reduces to the linear system: 
	\begin{align}
		\big(DF(z_*) - (\vec{k} \cdot \vec{\lambda}) I_n \big) W_{\vec{k}} &= -E_{\vec{k} }	 , \label{eq:CoHoLinear}
	\end{align}
	where $ W_{\vec{k}}, E_{\vec{k} } \in \C^n$ are the coefficients of the monomial $ \sigma^{\vec{k}} = \sigma_1^{k_1} \cdots \sigma_d^{k_d}$ in the Taylor expansions of $ W(\sigma)$ and $ E(\sigma)$ respectively. 
	For each $ |\vec{k}|\geq 2$ the coefficient $ E_{\vec{k}}$ depends only on the coefficients of $W$ of order strictly less than  $ |\vec{k}|$. 
	Hence the coefficients $ W_{\vec{k}}$ may be recursively solved for in \eqref{eq:CoHoLinear} so long as the matrix $(DF(z_*) -(\vec{k} \cdot \vec{\lambda}) I_n)$ is invertible. This fails to occur precisely when there is an internal  resonance, which is to say when: 
\begin{align*}
		\vec{k} \cdot \vec{\lambda} &= \lambda_i , & 1\leq i \leq n
\end{align*}
	for $ \vec{k} \in \N^d$, $|k|\geq 2$, and any eigenvalue $ \lambda_i $ of $DF(z_*)$. This can only occur if the eigenvalues $(\lambda_1  \dots  \lambda_n)$ are rationally dependent. 
	
%	 Thus, it can be said that dynamics on a generic invariant manifold are conjugate to linear flow. 
%	While the Hartman-Grobman theorem may guarantee that any hyperbolic equilibrium may be conjugated with its linearization. However this conjugacy need only be a homeomorphism, and these resonances are the essential obstruction to a higher regularity conjugacy. 

Returning to the dynamics resulting from the Galerkin truncation in \eqref{eq:Galerkin}, recall that the zero equilibrium has eigenvalues 
$\lambda_n =  -e^{i \theta} n^2 $ for integers $ 0 \leq n \leq N$. 
%As we increase the Galerkin truncation, we obtain more resonance in the eigenvalues. 
%Since the eigenvalues are given by $ n^2 $ for $n \in \Z$, 
Thus, resonances occur whenever we can write a square integer as a sum of other square integers:
\[
m^2 = \sum_{n=1}^{m-1} k_n n^2  
\]
for non-negative integers $k_n$ with $ k_1 + \dots+ k_{m-1} \geq 2$. 
This happens abundantly often!  For example $ m^2 = m^2 \cdot 1^2$ or $ m^2 = (m^2-n^2) \cdot 1^2 + 1 \cdot n^2$ for any $ 1 < n < m$. 
While resonances are an obstruction to conjugating the internal dynamics of $\cW$ to a linear system, a parameterization may still be obtained for more nontrivial internal dynamics.

 For an example,  consider the $N=3$ Galerkin truncation of \eqref{eq:Galerkin} given by:  
	\begin{align} \label{eq:Galerkin4}
%		e^{ -i \theta } \dot{a}_0 &= a_0^2 + 2 \sum_{k=1}^3 a_k^2\\
%		e^{ -i \theta } \dot{a}_1 &= -a_1  + 2 \sum_{k=1}^3 a_{k-1} a_{k} \\
		e^{ -i \theta } \dot{a}_0 &= a_0^2 + 2 a_1^2 + 2 a_2^2 + 2 a_3^2 \nonumber\\
		e^{ -i \theta } \dot{a}_1 &= -a_1  + 2 a_0 a_1 + 2 a_1 a_2 + 2 a_2 a_3 \\
		e^{ -i \theta } \dot{a}_2 &= -4a_2  + 2 a_0 a_2 +  a_1^2  + 2 a_1 a_3 \nonumber\\
		e^{ -i \theta } \dot{a}_3 &= -9a_3  + 2 a_0 a_3 + 2 a_1 a_2 \nonumber
	\end{align}
Again, we define $ \cW(0)$ as the unique invariant manifold tangent to  $\{0\} \times \C^3$. The equilibrium's nonzero eigenvalues are $\lambda = e^{ -i \theta }  (  -1,-4,-9) $, and we obtain the resonances:  
\begin{align}
	(4,0,0)\cdot\lambda&= \lambda_2 ,&
	(1,2,0)\cdot\lambda&= \lambda_3 ,&
	(5,1,0)\cdot\lambda&= \lambda_3 ,&
	(9,0,0)\cdot\lambda&= \lambda_3  \label{eq:Resonances}
\end{align}
 One can deal with  resonances by adjusting the function $f$, that one conjugates the dynamics to, see for example  \cite{van2016computing}. 
In our   code \cite{codesUnstableBlowupNLS} we calculate a parameterization  $W$ with rational coefficients up to order $20$ such that the dynamics are conjugate to:  
\begin{align}
	\mu^{-1} \dot{\sigma}_1  &=  - \sigma_1  \nonumber \\
	\mu^{-1} \dot{\sigma}_2  &=  -4 \sigma_2 + \tfrac{1}{3} \sigma_1^4 \label{eq:ConjugateSystem}\\
	\mu^{-1} \dot{\sigma}_3  &=  -9 \sigma_3 -   \sigma_1 \sigma_2^2  +   \tfrac{19}{24}   \sigma_1^5 \sigma_2 +\tfrac{11}{81} \sigma_1^9  \nonumber 
\end{align}
where $ \mu = e^{i \theta}$. The nonlinear terms in \eqref{eq:ConjugateSystem} are in direct correspondence with the resonances in \eqref{eq:Resonances}. 
The specific  coefficients for the higher order terms, such as $ \frac{1}{3} \sigma_1^4$ and $ \tfrac{19}{24}   \sigma_1^5 \sigma_2 $, are not able to be determined in advance and requires one to first solve the cohomological equation \eqref{eq:Cohomological} for all of the lower order terms.  
Furthermore, by using the parameterization method to solve for the coefficients up to an order past the point of all possible resonances, we ensure that the cohomological equation \eqref{eq:Cohomological} with the choice of $ f_k =0 $ can be solved for all $ k \geq 9$.

 While nonlinear, the internal dynamics of $\cW$  determined by \eqref{eq:ConjugateSystem} is  integrable, and for initial data $ \sigma(0) = (\gamma_1,\gamma_2,\gamma_3)$ its solution is explicitly given by:  
\begin{align}
	\sigma_1(t) &= \gamma_1 e^{-\mu t} , \nonumber  \\
	\sigma_2(t) &= \left( \gamma_2 +\frac{\mu t}{3} \gamma_1^4\right)e^{-4\mu t} ,\label{eq:ConjugacySolution}\\
	\sigma_3(t) &=    \left(
	\gamma_3 
	-\mu t \gamma_1  \gamma_2^2   
	+ 
	\left(
\frac{19 \mu t}{24} - 	\frac{ \mu^2 t^2}{3}
	\right) \gamma_1^5 \gamma_2 
	+ 
	\left(
		\frac{11 \mu t}{81}
	+ 	\frac{19 \mu^2 t^2}{144}
	-	\frac{ \mu^3 t^3}{27}
	\right) \gamma_1^9 \right)e^{-9\mu t}   \nonumber 
\end{align}
Note that if $ -\pi/2 < \theta < \pi/2$ then $ Re(-\mu) <0$ whereby  all of these solutions exponentially decay. 

 However if $ \theta = \pm \pi/2$ then $ \mu = \pm i$, whereby the solutions in \eqref{eq:ConjugacySolution} are oscillatory with secular drift exciting the higher modes. 
 Indeed, while $|\sigma_1(t)| $ stays bounded for all $ t$, the higher modes grow like  $|\sigma_2(t)| \sim t |\gamma_1|^4$ and $|\sigma_3(t)| \sim t^3 |\gamma_1|^9$. 
 Hence small initial data in the lowest modes requires a long time before the higher modes are excited to a comparable level.  
For example if we consider $ \gamma(0) = ( \epsilon,0,0)$ with $ \epsilon \ll 1$, then it takes $ t = \epsilon^{-2.5}$ time until $ |\sigma_3(t)| \gtrsim |\sigma_2(t)|$, and  just $ t = \epsilon^{-2.6\bar{6}}$ time until $ |\sigma_3(t)| \gtrsim |\sigma_1(t)|$. This ordering of the growth in the first few modes qualitatively  matches the behavior we observe in the PDE dynamics, as shown in Figure \ref{fig:BlowupPlot}(c).

\begin{figure}[b!]
	\begin{center}
		\includegraphics[width=1\textwidth]{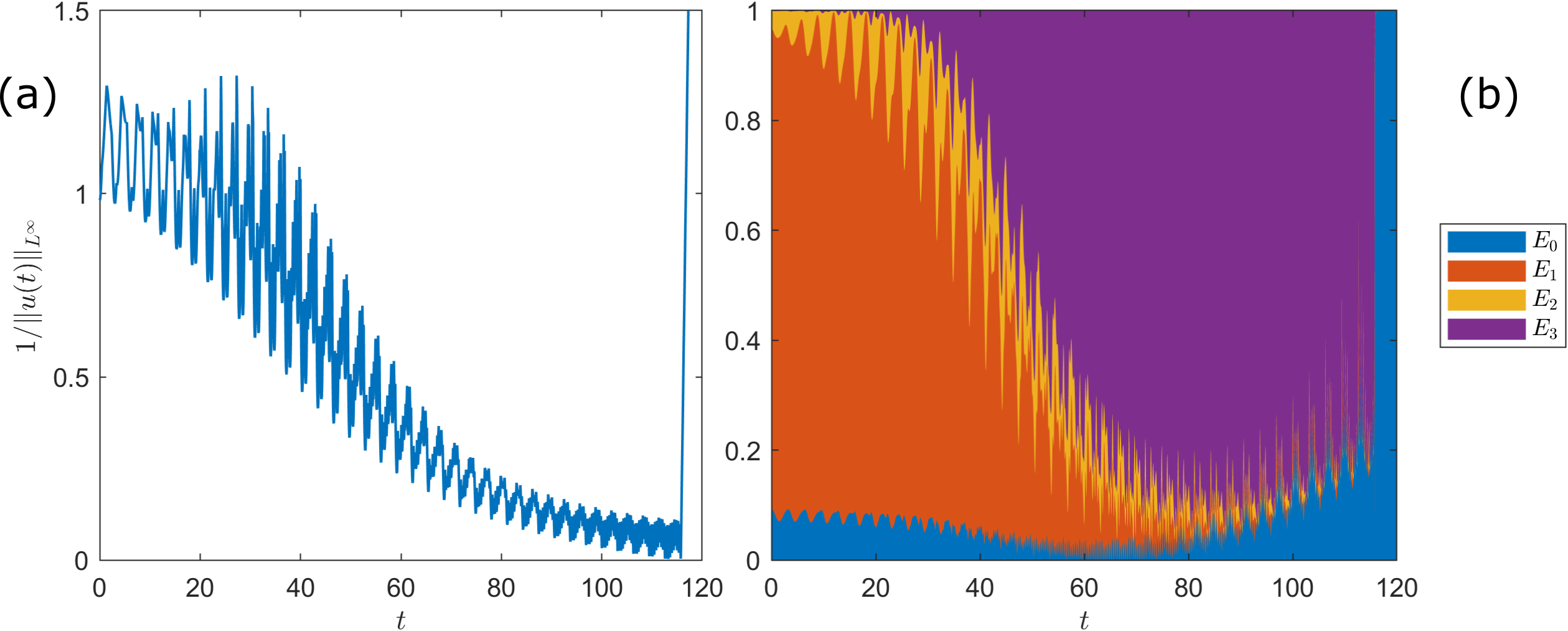} 
	\end{center}
	\caption{Solution to the four-mode Galerkin truncation in \eqref{eq:Galerkin4} with initial data \eqref{eq:FourModeInitialData} (a): The inverse norm $1/\| u^{(N)}(t)\|_{L^{\infty}}$ ; (b) Relative proportions $E_n$ of the Fourier modes, see \eqref{eq:ProporionalEnergy}.}
	\label{fig:ToySystemSimulation}
\end{figure}

Like with any power series, the parameterization of $\cW$ is only valid on the power series' radius of convergence.  
While the explicit solutions given in \eqref{eq:ConjugacySolution} never blow up, it may be possible for a solution on $ \cW$ to blowup after it leaves the local coordinate chart. 
Based on our computation, the radius of convergence seems to be about $0.80$.  
To investigate this, we consider the following initial condition selected on the invariant manifold $\cW$:  
\begin{align}
	\left(
	\begin{matrix}
		\sigma _1  \\ 
		\sigma _2  \\ 
		\sigma _3 
	\end{matrix}
	\right)
	&=
		\left(
	\begin{matrix}
		 \;0.4300654917290795 \\ 
		 -0.07398732057014827 \\ 
		\;\;0.00530826265454094
	\end{matrix}
	\right)
	,
	&W(\sigma )
	&=
	\left(
	\begin{matrix}
		-0.22301409257004942 \\ 
		0.5 \\ 
		0\\
		0
	\end{matrix}
	\right)
	\label{eq:FourModeInitialData}
\end{align} 
This initial condition was selected such that $ a_2(0)=a_3(0)=0$ as to mimic the initial data in \eqref{eq:InitialDataFamily}.

When integrated under \eqref{eq:Galerkin4} with $ \theta = \pi/2$ this trajectory undergoes a large growth in norm,  see Figure \ref{fig:ToySystemSimulation}(a). 
Like with Figure \ref{fig:BlowupPlot},  there is a steady  cascade pumping energy into the higher modes and the solution eventually attains a maximum value of $\|u(t)\|_{L^\infty} = 335$ at time $t = 113.3$.
However this solution does not appear to completely blowup. As can be seen in Figure \ref{fig:ToySystemSimulation}(b), the zero mode initially decreases relative to the other components until time $ t \approx 60$. After this point the zero mode grows to be more and more dominant.  
Eventually %t=115.73
the trajectory is attracted to invariant  manifold $\cW^0 = \{ \{a_i\} \in \C^4 : a_i=0 \forall i\neq0\}$ of spatially homogeneous solutions, after which it decays algebraically to the zero equilibrium, cf Theorem \ref{thm:TrappingRegion}.

Such a non-blowup could still be seen to be consistent with the existence of an unstable blowup. 
Our conjecture that blowup occurs on a codimension-1 manifold means that the initial conditions need to be precisely selected in order to be observed through direct numerical simulation. 
The accumulation of numerical error over a long interval of time may cause the trajectory to drift away from the set of data which blows up. 
Furthermore our initial condition cannot be guaranteed to be exactly on the manifold $\cW$. 
Given our 20th order Taylor approximation  with an assumed radius of convergence of $0.8$, we estimate the error of the initial point to be on the order of $10^{-5}$. 
While one could get a more precise initial point on the manifold by choosing a smaller initial condition $\sigma$, this would conversely increase the integration time and the associated errors.

\section{Apparent self-similar blowup profiles}
\label{sec:SelfSimilar}

The finite Fourier truncation model offers a heuristic explanation for why initial data on the stable manifold of the heat equation \eqref{eq:Heat} will grow larger, oscillate, and exhibit a cascade to the higher modes. 
However this analysis is localized at the zero equilibrium. While it is suggestive of how small initial data grows to become finitely large in norm, it does not explain how large initial data may grow without bound and blowup.

\begin{figure}[h]
	\begin{center}
		\includegraphics[width=1\textwidth]{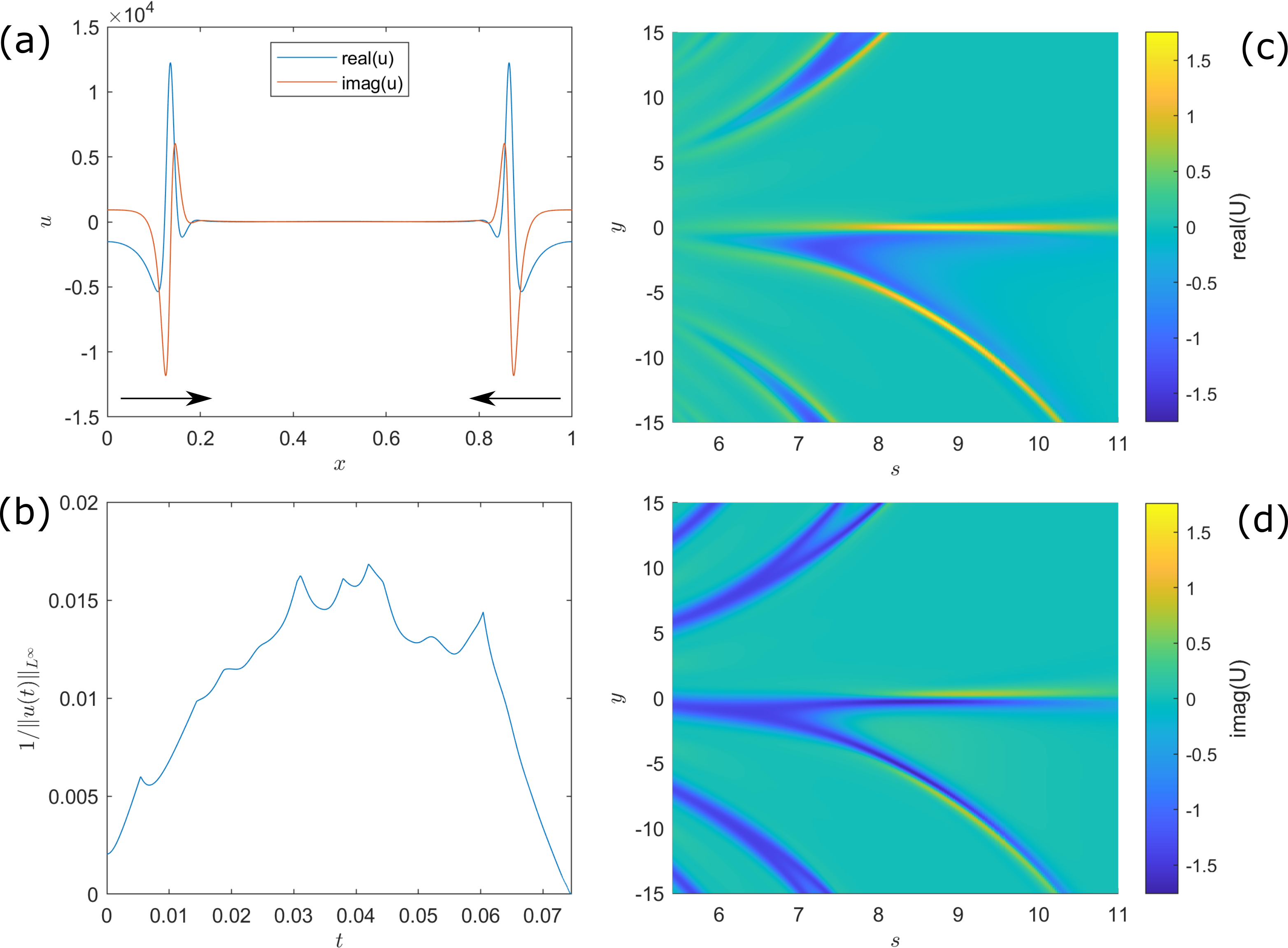} 
	\end{center}
	\caption{Numerical solution of \eqref{eq:NLS} using initial data from \eqref{eq:InitialDataA300}. (a) Real and imaginary components of $ u(t,x)$ when $t=0.0743$ ($s=8.95$%		8.947976107508760
		), depicted with the direction of movement of the blowup point(s); (b) 
		Inverse norm of solution $1/\| u(t)\|_{L^{\infty}}$;
		 (c-d) Real and imaginary components of  $U$ from \eqref{eq:SelfSimilarCoordsU} using self-similar coordinates with scaling  $(\alpha,\beta)=(1,1/2)$.	}
	\label{fig:Combined_A300}
\end{figure}

\begin{figure}[h]
	\begin{center}
		\includegraphics[width=1\textwidth]{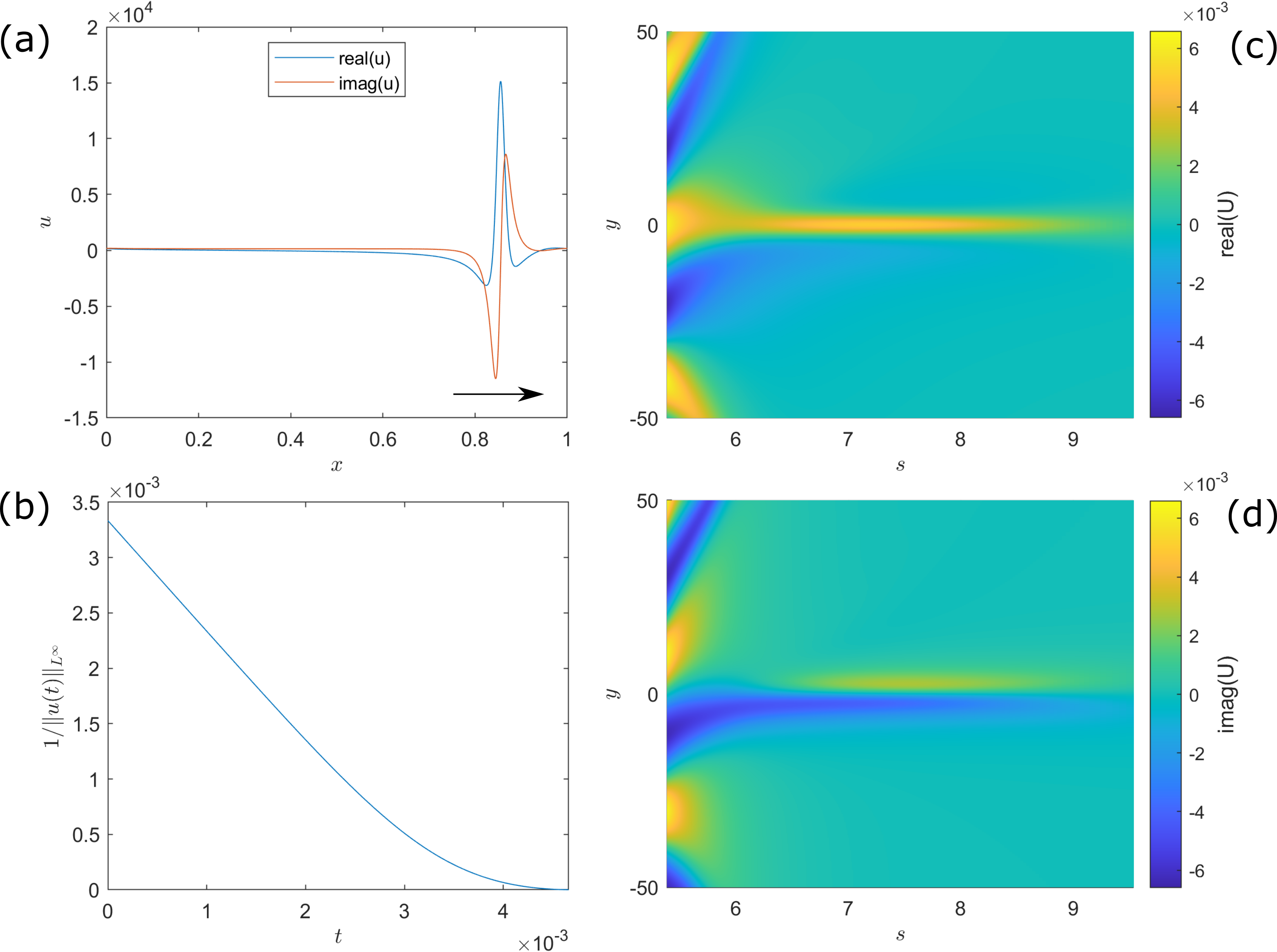} 
	\end{center}
	\caption{Numerical solution of \eqref{eq:NLS} using initial data from \eqref{eq:InitialDataE300}.  (a) Real and imaginary components of $ u(t,x)$ when $t=0.004$ ($s=7.47$%		7.466371566584477
), depicted with the direction of movement of the blowup point; (b) 
		Inverse norm of solution $1/\| u(t)\|_{L^{\infty}}$;
		 (c-d) Real and imaginary components of   $U$ from \eqref{eq:SelfSimilarCoordsU} using self-similar coordinates with  scaling  $(\alpha,\beta)=(2,1)$. }
	\label{fig:Combined_E300}
\end{figure}

%While the data $ u_{30}$   grows very large the late stage behavior near blowup is immensely complicated. 
To focus on the dynamics of blowup  we consider larger initial data. With  consideration to Conjecture \ref{conj:ParabolicManifold}, we fix the following real initial data: 
\begin{align} \label{eq:InitialDataA300}
	u_{300}(x) &= 300 \cos(2 \pi x) + A_{300}, & A_{300} &= -189.286840601635,
\end{align}
%More precisely, using a bisection method, we selected $ A_{300}$ such that $u_{300}$ would be on the stable manifold $\cW(0)$ for the nonlinear heat equation \eqref{eq:Heat}.   
see Figure \ref{fig:Combined_A300}. 
We used  computational parameters of 4096 Fourier mode truncation and a time step of $h=10^{-7}$. 
For comparison we also consider monochromatic initial data: 
\begin{align}\label{eq:InitialDataE300}
	u_{300}^{mc}(x) &= 300 e^{2 \pi i x},
\end{align}
see Figure \ref{fig:Combined_E300}. 
The initial data in \eqref{eq:InitialDataE300} is guaranteed to blowup due to the Theorem \ref{prop:Intro_Blowup}, (noting also that $300 > 6 (2\pi)^2 \approx 236.9$). 

%Note that the positive Fourier mode assumption is somewhat removed from the real initial data case, as there one has $ a_k = \bar{a}_{-k}$.  

As predicted by   Conjecture \ref{conj:ParabolicManifold} and Theorem \ref{prop:Intro_Blowup}, the numerical solutions of the initial data \eqref{eq:InitialDataA300} and \eqref{eq:InitialDataE300} appear to blowup.  
Plotted in Figures \ref{fig:Combined_A300}-\ref{fig:Combined_E300}(a) are the numerical solutions after a certain amount of time.
The monochromatic initial data yields a single blowup point moving to the right. 
However the real initial data produces a solution which at first appears to blowup at the stationary point $x=0$ (or the periodically equivalent point $ x=1$). 
However at around time $ t=0.0735$ %($ s = 0.07$)
 this blowup point appears to bifurcate into two non-stationary blowup points traveling in opposite directions.

 Moreover,  the solutions appear to undergo  self-similar blowup.  
 In a general context, this occurs when a  blowup   solution $u$ to a given PDE can  be regularized using  a self-similar change of variables into a new solution $U$  with   nontrivial limiting behavior \cite{eggers2008role,quittner2019superlinear}.  
For a given blowup time $T$ and pair of scaling parameters $ (\alpha,\beta)$  define: 
\begin{align}
u(t,x) &= \frac{1}{(T -t)^\alpha} U(s,y)  \label{eq:SelfSimilarCoordsU}
%,&
%	y&= \frac{x-\xi(t)}{(T -t)^\beta} ,  \;\;
%	s= - \log( T - t)
\end{align}
using the self-similar coordinates: 
\begin{align}
	y&= \frac{x-\xi(t)}{(T -t)^\beta} 
	&
	s&= - \log( T - t)
\end{align}
where $ \xi(t)$ denotes the moving location of the spatial blowup point. 
%Numerically, it appears that $  \xi'(t)  $ limits to a constant value as $ t \to T$.
Note that if $ \alpha >0$ and $ \limsup_{s \to \infty } \| U(s) \|_{L^\infty} > 0$,  then $u(t)$ blows up as $ t \to T$.    While the dynamics of $ u$ are singular as $ t \to T$, ideally the dynamics of $U$ are  well behaved as $ s \to \infty$.

 In Figures \ref{fig:Combined_A300}-\ref{fig:Combined_E300}(b) the inverse norm $ 1 / \| u(t)\|_{L^{\infty}}$  can be seen to approach zero (thus blowing up). 
 The solution with real initial data first decreases in norm and oscillates to a small extent before it begins a  path towards blowup. 
 In contrast, the solution with monochromatic initial data immediately increases in norm and accelerates in the later stage.  
  
 To estimate the scaling parameter $ \alpha$, we fit the data in Figures \ref{fig:Combined_A300}-\ref{fig:Combined_E300}(b) to the equation $ y(t) = C_0 |T-t|^\alpha$.  
 For the solution with real initial data, fitting the data over the time interval $[0.070,0.074]$ yielded $\alpha =1.1457%\pm0.0024
 $ with an $R^2$ value of $0.9999$. 
  For the solution with monochromatic initial data, fitting the data over the time interval $[0.038,0.045]$ yielded $\alpha =2.0098%\pm0.0638
  $ with an $R^2$ value of $0.9992$. 
  Note that even though the two blowup profiles in Figures \ref{fig:Combined_A300}-\ref{fig:Combined_E300}(a) look qualitatively similar, they appear to obey distinct rates of blowup. 
Additionally, while the $R^2$ values are  suggestive of a good fit, the residual errors are not normally distributed,  indicative of systematic bias. Practically speaking, we observe that fit parameters are quite  sensitive to the time range over which the data is fit.

To investigate  whether these numerical solutions obey self-similar scaling, we perform a  change of variables in \eqref{eq:NLS} using the self-similar coordinates, resulting in the following equations that govern the self-similar dynamics: 
\begin{align} \label{eq:SelfSimilar}
		i 
\left(  \alpha U + \partial_s U +  \beta y  \partial_y U  \right)
&=  e^{(2 \beta -1)s} \partial_{yy} U
 + 
e^{(\beta-1)s}  i\dot{\xi}
\partial_y U
+
e^{(\alpha-1)s} 
 U^2 
\end{align} 
For consistency in scaling, balance laws suggest that $ \alpha =2 \beta$.   
Rounding to the closest integer our statistical fits for $\alpha$, we use scaling parameters  $ (\alpha,\beta)=(1,1/2) $ for the real initial data and  $ (\alpha,\beta)=(2,1)$ for the monochromatic initial data. 
Making the change of variables into these self-similar coordinates, we plot the solutions $U(s,y)$ in Figure  \ref{fig:Combined_A300}-\ref{fig:Combined_E300}(c-d)  for the real and imaginary components of the  solutions. 

%While the profiles appear to converge to a profile decently well,  the amplitude of the profile fades for large $s$. 

In Figure  \ref{fig:Combined_A300}-\ref{fig:Combined_E300}(c-d) the self-similar solutions $U(s,y)$ appear to be decently regularized. For smaller values of $ s$, periodic copies of the blowup profile may be observed to diverge away from $y=0$, as would be expected from a PDE posed on a periodic domain. 
For larger values of $s$ however the blowup profile appears to fade, a likely result of an imprecise selection of the blowup time, and the significant numerical error which accumulate as the blowup time is approached. 
While the solution with real initial data appears to grow according to a power law  starting at $ t = 0.06$ ($s=4.25$), the blowup profile qualitatively changes. 
This may be most prominently seen in Figure \ref{fig:Combined_A300}(d) when comparing the imaginary component of $U(s,y)$ over the regions $ s \in [6,7]$ and $ s \in [8,10]$, and is suggestive of non-trivial self-similar dynamics. 
 
Alternatively,  the simplest blowup scenario would be a so-called Type-I blowup \cite{eggers2008role,quittner2019superlinear}, whereby $U$  converges to an equilibrium with scaling parameters $(\alpha,\beta)=(1,1/2)$ and \eqref{eq:SelfSimilar} simplifies to:  
	\begin{align}	\label{eq:TypeI}
		i 
		\left(    U + \partial_s U +    \tfrac{y}{2}  \partial_y U  \right)
		&=   \partial_{yy} U
		+ 
		e^{-s/2}  i\dot{\xi}
		\partial_y U
		+ 
		U^2 
	\end{align}  
In this case, as $ \lim_{s \to \infty} e^{-s/2} =0$,   an equilibrium solution would satisfy the following second order ODE:
\begin{align}	0
	&=   \partial_{yy} U - \tfrac{i}{2} y  \partial_y U -iU
	+ 
	U^2,  & \lim_{y \to \pm \infty} U(y) \in \{0,i\}
	\label{eq:BlowupProfile2ndOrder}
\end{align}
Despite our efforts we have not found a nontrivial solution to \eqref{eq:BlowupProfile2ndOrder}. 
We suspect that some other self-similar scaling ansatz is needed, such as a logarithmic correction \cite{eggers2008role,bricmont1994universality,merle1997stability}.

 \section{Conclusion}

We have presented numerical evidence of real initial data to \eqref{eq:NLS} which blows up in finite time, in contradiction with Conjecture  \ref{conj:GWP}.  
By tracking solutions to the 1-parameter family of initial data in \eqref{eq:InitialDataFamily}, we were able to identify initial conditions leading to blowup. 
While tracking the $L^\infty$ norm of solutions in Figure \ref{fig:ContinuationInZeroAverage} proved sufficient to identify blowup solutions, measuring when solutions entered a trapping region of zero provided a much more robust method of identification. 
 
% that as stated is false. 
% It still may be such that that the conjecture is true in a generic sense (almost all real initial data is GWP).  
% Through Theorem \ref{thm:TrappingRegion} we are able to restrict this search for a given initial data to searching through its average value $A$. 
% By tracking the $L^\infty$ norm of solutions as in Figure \ref{fig:ContinuationInZeroAverage} one can identify using a fine enough search initial conditions which appear to blowup. 
% In contrast we were able to more robustly identify blowup by tracking how long a given initial takes to enter the trapping region of the zero equilibrium. 
 
The blowup solution identified in Section \ref{sec:EvidenceOfBlowup} has two important qualitative features: (i) on shorter time scales, the solution  periodic oscillates with progressive excitement of the higher modes; (ii) on longer time scales the solution steadily grows, eventually leading to blowup.  
In Section \ref{sec:Resonances} we provide a heuristic explanation for the mechanisms behind feature (i). 
Namely, we use the parameterization method to analyze a submanifold $\cW$ of the center manifold of the zero equilibrium. 
While solutions on $\cW$ do oscillate, a secular drift due to a resonance of eigenvalues induces a forward cascade, whereby the lower modes excite the higher modes. 

This analysis leads us to propose Conjecture \ref{conj:ParabolicManifold}, that real initial data to \eqref{eq:NLS} will blowup if and only if it lies on the strong-stable manifold of $0$ for the nonlinear heat equation \eqref{eq:Heat}.  
We note that the finite Galerkin truncation model only has a finite number of eigenvalue resonances, and is thus amenable to the parameterization method. %, however the full PDE has an infinite number of resonant eigenvalues. 
While empirical evidence supports Conjecture \ref{conj:ParabolicManifold}, we cannot hope to prove it using the parameterization method due to the infinite number of resonances  in the full PDE.

To analyze the later stage of blowup in \eqref{eq:NLS} we performed a self-similar analysis in Section \ref{sec:SelfSimilar}, comparing blowup solutions starting from both real and monochromatic initial data. 
While both solutions appear to exhibit self-similar blowup, key questions remain. 
For example, what is the exact scaling rate for the solution starting with real initial data? 
And does the self-similar solution limit to an equilibrium or some more complicated dynamical object? 
 
\section*{Acknowledgments}

The author would like to thank Panayotis Kevrekidis and Javier G\'omez-Serrano for   fruitful discussions regarding this work.

\bibliography{Bib_NLS}
\bibliographystyle{alpha}

\end{document}